\documentclass[]{aiaa-tc}

\usepackage{mathpazo}
\usepackage{setspace}
\usepackage{amsmath}
\allowdisplaybreaks
\usepackage{graphicx}
\usepackage{booktabs}

 \title{Routing Multiple Unmanned Vehicles in GPS-Denied Environments}

 \author{ \normalsize
   Bingyu Wang%
   \thanks{PhD Student, Dept. of Mechanical Engg., Texas A\&M University, College Station, TX, USA.}
   \;
   Sohum Misra
   \thanks{PhD Student, Dept. of Aerospace Engg., University of Cincinnati, Cincinnati, OH, USA.}
   \;
    Sivakumar Rathinam
  \thanks{Associate Professor, Dept. of Mechanical Engg., Texas A\&M University, College Station, TX, USA.}
  \;
  Rajnikant Sharma
  \thanks{Assistant Professor, Dept. of Aerospace Engg., University of Cincinnati, Cincinnati, OH, USA.}
  \;
  Kaarthik Sundar%
    \thanks{Postdoctoral Researcher, Center for Nonlinear Studies, Los Alamos National Laboratory, Los Alamos, NM, USA}
 }

\usepackage{amsmath}
\usepackage{amssymb}
\usepackage{graphicx}

\begin{document}

\maketitle

\begin{abstract}
This article aims to develop novel path planning algorithms required to deploy multiple unmanned vehicles in Global Positioning System (GPS) denied environments. Unmanned vehicles (ground or aerial) are ideal platforms for executing monitoring and data gathering tasks in civil infrastructure management, agriculture, public safety, law enforcement, disaster relief and transportation. Significant advancement in the area of path planning for unmanned vehicles over the last decade has resulted in a suite of algorithms that can handle heterogeneity, motion and other on-board resource constraints for these vehicles. However, most of these routing and path planning algorithms rely on the availability of the GPS information. Unintentional and intentional interference and design errors can cause GPS service outages, which in turn, can crucially affect all the systems that depend on GPS information. This article addresses a multiple vehicle path planning problem that arises while deploying a team of unmanned vehicles for monitoring applications in GPS-denied environments and presents a mathematical formulation and algorithms for solving the problem. Simulation results are also presented to corroborate the performance of the proposed algorithms.
\end{abstract}

\section*{Nomenclature}

\begin{tabbing}
  XXX \= \kill
  $I$ \> set of depots \\
  $J$ \> set of targets \\
  $LM$ \> set of potential landmark locations \\
  $G$ \> graph $G$ with vertex set $V$ and edge set $E$ \\
  $S$ \> any vertex subset \\
  $F$ \> any edge subset \\
  $v$ \> speed of a vehicle \\
  $\psi$ \> heading of the vehicle \\
  $ \omega $ \> yaw rate \\  
 \end{tabbing}

\section{Introduction}

\subsection{Motivation and related work}

This article considers a routing problem involving a team of Unmanned Vehicles (UVs) in GPS-denied environments. UVs have been commonly used in several civil and military applications over the last decade. A suite of path planning algorithms that can handle heterogeneity, motion and other onboard resource constraints for UVs have been developed previously in the literature \cite{sundar2016generalized,sundar2016formulations,toth2014vehicle}. However, most of these algorithms rely on the availability of the GPS information. This makes UVs vulnerable to GPS jamming and spoofing \cite{kerns2014unmanned,hoey2005civil,carroll2003vulnerability}; moreover, most indoor environments and many parts of the urban canyon do not have access to GPS and even if available, the access is intermittent and not reliable. Hence, routing in a GPS-denied or GPS-restricted environment is an important area of research \cite{manyam2016gps,sundar3}.

In the absence of GPS information, on-board range or bearing angle sensors such as LIDAR and camera can be used to infer relative position measurements of a vehicle to known landmarks (LMs). The process of estimating the relative position of a vehicle with respect to the sensed landmarks is known as localization \cite{sharma2013bearing}. The process of both building a map of the environment and estimating the relative position of a vehicle in the map is known as Simultaneous Localization And Mapping (SLAM). Typically in a SLAM problem, the path followed by the vehicle is mainly focused towards building a good map or localizing the relative position of the vehicle with respect to the map. In this article, we consider the problem of planning good paths for the vehicles and placing the landmarks suitably so that the vehicles can visit the given set of waypoints (also referred to as targets) as quickly as possible in addition to satisfying some localization constraints as the vehicles traverses their paths. We refer to this problem as a Multiple Vehicle path planning Problem with Localization Constraints (MVPLC). 

In the presence of GPS information, given a set of targets and a single vehicle, the problem of finding a shortest path for visiting each target in a given set of targets is called the Traveling Salesman Problem (TSP) and is computationally hard to solve. MVPLC can be viewed as a generalization of the TSP where additional constraints are imposed on the choice of the vehicle paths. In addition, MVPLC also aims to place the landmarks suitably so that the vehicles can use the landmarks efficiently for localization. The placement of LMs given a set of paths for multiple vehicles was first considered in Rathinam et. al \cite{rathinam2015multiple}. The authors formulated the landmark placement problem as a multiple vehicle path covering problem and presented an approximation algorithm using geometric arguments to address the problem. In the most recent work, \cite{sundar3} the authors offered fast heuristics to the landmark placement problem, and found an optimal solution to the single vehicle routing problem. This article is a natural extension of the existing work \cite{sundar3} to multiple vehicles.

\subsection{Problem statement}
Given a fleet of homogeneous UVs starting at distinct depots, a set of target locations, and a set of potential locations where the LMs can be placed, the MVPLC aims to find a route for each of the UVs and a subset of potential LM locations where the LMs are placed, such that the following conditions are satisfied:
\begin{enumerate}
\item the route assigned to each UV ends at its respective starting location (depot) of the vehicle,
\item each target is visited at least once by any of the UVs,
\item each UV should be able to estimate its position and orientation from the LMs as it traverses its route, and
\item the sum of the total travel cost and the placement cost of LMs used is a minimum.
\end{enumerate}

\section{Mathematical Formulation} \label{sec:form}

We formulate MVPLC as a mixed-integer linear program (MILP). The notation and some basic results that will be used throughout the paper are introduced now. We define MVPLC on a set of clients \textit{J}, a set of potential depots \textit{I}, and a set of potential LM locations \textit{LM}. Denote $ |J|=q, |I|=p $, and $ |LM|=lm $ and assume that they are positive integers. Let $G = (V ,E)$ be an undirected graph where $V = I \cup J$, and $ E = {(i, j ) : \forall i \in V,\forall j \in J } $ (note that \textit{E} does not include any edge between depots). The cost of edge $ e = (i, j ) $ is denoted by $ c_{ij}=c_e $. To make use of methodologies in solving standard TSPs, assume the costs of all edges satisfy the triangular inequality, i.e. $c_{ij}+c_{jk}>c_{ik}, \forall i,j,k \in V $.

For each edge $ e = (i, j ), i, j \in J $, we define one binary variable $x_{ij}$ which takes the value 1 if the edge \textit{e} is traveled by one route and 0 otherwise. For each edge $e = (i, j ), i \in I, j \in J$ we define a variable $x_{ij}$ which takes the value 2 if one vehicle does a trip between depot \textit{i} to client \textit{j} and immediately comes back to the depot (this is called a return trip), the value 1 if the edge \textit{e} is traveled once by one vehicle, and 0 otherwise. For two vertex subsets $ S,S'\subseteq V $, define $ (S:S')=\lbrace(i, j):i \in S, j \in S' \rbrace $ . Given a vertex subset, $ S\subseteq V $, denote $ \delta(S)=(S:V\setminus S) $ and $ \gamma(S)=\lbrace(i, j)\in E: i,j \in S \rbrace $. If S is a one-vertex singleton, use $ \delta(v)$ instead of $ \delta (\lbrace v \rbrace)$. Finally, for $F \subseteq E$, define $x(F) =\sum_{(i, j)\in F} x_{ij}$. 

For each potential LM location $ k \in LM $, we define one binary variable $d_k$ which takes the value 1 if the location \textit{k} is used to place an LM and 0 otherwise. For each edge $ e\in E $, let $L_e$ denote a subset of feasible locations to place LMs, such that when a vehicle is travelling via $e$, it can be localized with the help of LMs in $L_e$. We assume that the cost of placing an LM is unitary; but even if it's not, the problem structure remains exactly the same.

The formulation of the MVPLC is as follows:
\begin{flalign}
& \min \quad \sum_{e\in E} c_e x_e + \sum_{k=1}^{lm} d_k \text{ subject to: } \notag \\
& x(\delta(j))=2, \quad \forall j\in J, \label{e:1} \\
& x(\delta(S))\geqslant 2, \quad \forall S\subset J, \label{e:2} \\
& \sum_{i\in I'} x_{ij} +2x(\gamma(S\cup \lbrace j,l\rbrace)) + \sum_{k\in I\setminus I'} x_{kl} \leqslant 2|S|+3 \quad \forall j,l\in J, S \subseteq J\setminus \lbrace j,l\rbrace, S \neq \varnothing; I'\subset I, \label{e:3} \\
& \sum_{i\in I'} x_{ij} +3x_{jl} + \sum_{k\in I\setminus I'} x_{kl} \leqslant 4, \quad \forall j,l\in J, I'\subset I,  \label{e:4} \\
& \sum_ {k\in L_e} d_k \geqslant 2 x_e, \quad \forall e\in E, \label{e:5} \\
& x_{ij}\in \lbrace 0,1,2\rbrace , \quad \forall i\in I, \forall j\in J, \label{e:6} \\
& x_{ij}\in \lbrace 0,1\rbrace, \quad  \forall i\in J, \forall j\in J, \label{e:7} \\
& d_k \in \lbrace 0,1\rbrace , \quad \forall k\in LM. \label{e:8}
\end{flalign}

Degree constraints in \eqref{e:1} ensure that all the targets are visited exactly once by the set of routes. Inequalities \eqref{e:2} are the very well-known subtour elimination inequalities; they can also be written in a complementary fashion. Inequalities \eqref{e:3}--\eqref{e:4}, called path elimination constraints, were first introduced in Laporte et al.\cite{laporte1986exact}, and prevent solutions that include a path starting at one depot and ending at another depot. The validity of them has been discussed in \cite{benavent2013multi}. Inequalities \eqref{e:5}, the localization constraints, state that there must be at least two landmarks in the proximity of any edge traveled by the vehicles. The vehicles essentially will use these nearby landmarks to localize its position as it travels along the edge. These constraints stem from the results on estimation. \cite{sharma2013bearing}\\

\section{Solution Methodology}

We will now outline our implementation methodology which will be used to corroborate the performance of the proposed algorithms in both hardware in the loop simulations and field experiments. The landmarks will be placed at suitable locations, and the paths will be assigned to the vehicles using the solution obtained by solving the mathematical formulation. The sequence in which the vehicles visits the targets provides the way points for the paths. We use exteroceptive sensors on vehicles to provide the bearing information, from the LMs placed on the locations chosen by the optimization problem.

We now present the kinematic model of the vehicle in state-space form. We make the following asumptions:
\begin{enumerate}
\item the vehicle travels at a constant altitude with a constant velocity;
\item the vehicle cannot move backwards;
\item the sensing range of the vehicle's exteroceptive sensor, denoted by $\rho_s$, is constant;
\item the exteroceptive sensors equipped on the vehicle have a complete $360^{\circ}$ field of view, to measure the kinematics of the vehicle with respect to the known LMs.
\end{enumerate}
\begin{flalign*}
\dot{\mathbf{x}} = f(\mathbf{x,u})\triangleq 
\begin{bmatrix}
v\cos(\psi)\\
v\sin(\psi)\\
\omega
\end{bmatrix}
\end{flalign*}
where $ v $ is the magnitude of velocity of the vehicle, $ \psi $  is the orientation (heading angle) of the vehicle and $ \omega $ is the rate of change of heading with respect to time. The vector $\mathbf{x}$ is the vector of state variables given by $(x,y,\psi)^{T}$; here x and y denote the position of the vehicle and $ \psi $  gives the heading information. The control input vector $\mathbf{u}$ for the vehicle consists of $ v $ and $ \omega $, which are detected by the interoceptive sensors.

In this article, the Extended Information Filter (EIF), which is the information form of the Extended Kalman Filter (EKF), is used to estimate the states of the system $\mathbf{x}$, as the vehicles travel along their assigned paths. In our previous work \cite{sharma2013bearing}, we have shown that using the Hermann-Krener criterion for local observability, if each edge traversed by the vehicle is covered by at least two installed LMs (\textit{i.e.} two LMs within the sensing range of the vehicle), the system will be observable. We also remark that, having edges associated with more than two LMs provides more information to the vehicles, thereby quickening the convergence of the estimation algorithms. This method of state estimation is referred to as \textit{bearing-only localization}.

\begin{figure}
\centering
\includegraphics[width=100mm,scale=0.5]{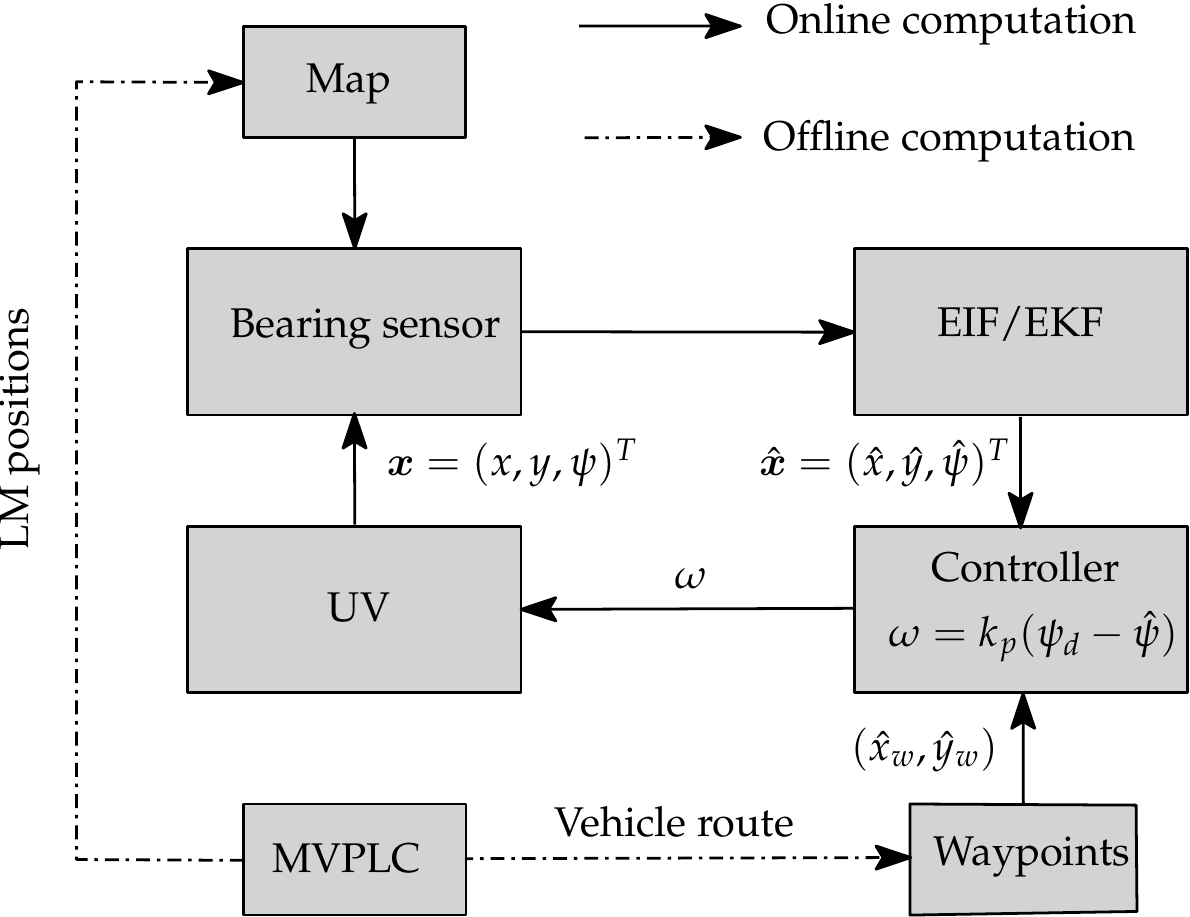}
\caption{The block diagram of the system architecture used for state estimation.}
\label{fig:arch}
\end{figure}

The states are provided as input to a proportional controller which computes the corresponding turn rate ($ \omega $, control input) to the vehicle. Taking from our previous work \cite{sundar3}, we will choose the controller gain such that a minimum distance requirement to each way point is satisfied. This in turn implies that for the vehicle to switch from the current way point to the next one, the vehicle needs to come in close proximity to the current way point (within a specified minimum distance) for the switching to take place. We do assume though that the vehicle can have a very high turn rate. The overall proposed architecture is shown in Figure \ref{fig:arch}. The dashed lines in the architecture denote that the computation is performed offline. 

\section{Algorithm Overview}

\subsection{Branch-and-cut Algorithm} \label{sec:bandc}

For the sake of completeness, a detailed pseudo-code of the standard algorithm with respect to the MVPLC is given. Let $ \bar{\tau} $ denote the optimal solution to an instance of MVPLC.\\

\noindent S\textsc{tep} 1 (Initialization). Set the iteration count $ t \leftarrow 1 $ and the initial feasible solution objective $ \bar{\alpha} \leftarrow +\infty $. The initial linear sub-problem is then defined by formulation in Sec. \ref{sec:form} without the sub-tour elimination constraints in Eq. \eqref{e:2}, the path elimination constraints in \eqref{e:3}--\eqref{e:4}, and the binary restrictions on the variables relaxed. The initial sub-problem is solved and inserted in a list $L$.\\

\noindent S\textsc{tep} 2 (Termination check and sub-problem selection). If the list $ L $ is empty, then stop. Otherwise, select a sub-problem from the list $ L $ with the smallest objective value.\\

\noindent S\textsc{tep} 3 (Sub-problem solution). $t\leftarrow t+1 $. Let $ \alpha $ denote the objective value of the sub-problem solution. If $ \alpha \geqslant \bar{\alpha} $, then proceed to S\textsc{tep} 2. If the solution is feasible for the MVPLC, then set $ \bar{\alpha} \leftarrow \alpha $, update $ \bar{\tau} $ and proceed to S\textsc{tep} 2.\\

\noindent S\textsc{tep} 4 (Constraint separation and generation). Using the separation algorithm for the sub-tour and path elimination constraints, identify the violated constraints in Eq. \eqref{e:2}--\eqref{e:4}. Add the violated constraints to the initial linear sub-problem and proceed to S\textsc{tep} 3. If no constraints are generated, then proceed to S\textsc{tep} 5.\\

\noindent S\textsc{tep} 5 (Branching). Create two sub-problems by branching on a fractional $ x_e $ or $ d_k $ variable. Then insert both the sub-problems in the list $ L $ and go to S\textsc{tep} 2. \\

The following paragraphs detail the algorithms to identify violated sub-tour elimiation and path elimination constraints in S\textsc{tep} 4 of the algorithm. 

\subsection{Separation of sub-tour elimination constraints}

Let $G^{*} = (V^{*},E^{*})$ denote the support graph associated with a given fractional solution $(x^{*}, d^{*})$ i.e., $V^{*}  = J$ and $E^{*} := \lbrace e \in E : x_e^{*}> 0 \rbrace $. Here, $x$ and $d$ are the vector of decision variable values in MVPLC. Next, we examine the connected components in $G^{*}$. Each connected component that does not contain all the targets in $J$ generates a violated sub-tour elimination constraint for $S = C$. If the number of connected components is one, then the most violated constraint of the form $ x(\delta(S))\geqslant 2 $ can be obtained by computing the global minimum on a capacitated undirected graph $G^{*}$; let the cut be denoted by $(S, V^{*} \setminus S)$. $S$ defines a violated sub-tour elimination constraint if the value of the cut is strictly less than $2$. This algorithm has been proposed as early as in 1990 \cite{padberg1990facet} for the traveling salesman problem.

\subsection{Separation of path elimination constraints}

According to the algorithm described by \cite{labbe2004branch}, we can keep using the idea of sequentially shrinking edges, whose weight is equal to one in the support graph. In our problem, only edges with both end-points in $J$ were shrunk. If the shrunk graph contains a vertex, say $v$, linked to more than one depot, we check whether a path elimination constraint is violated. Note that vertex $v$ will in fact correspond to a path in $G(\bar{x})$ where all the edges have weight $1$. Then identifying the extremes of this path as the clients $j$ and $l$, is quite easy to check if there is a subset of depots $I'$ for which the corresponding path elimination inequality is violated by $ \bar{x} $. Note that if $ \bar{x} $ is integer this procedure will find a violated path constraint if any exists.



\section{Computational and Simulation Results}

This section presents extensive computational and simulation results for all the algorithms developed thus far. All the computational experiments were performed on an ASUS X556U laptop with a 2.9 GHz Intel Core i7 processor and 8 GB RAM using CPLEX 12.7.1 as a mixed-integer linear programming solver.

\subsection{Instance generation}

We simulate cases for two, three, and four depots. The number of targets (including the depots), $ |V| $, for all the test instances was chosen from the set $ \lbrace $20, 25, 30$ \rbrace $. For each value of $ |V| $, $20$ random instances were generated. The targets were randomly placed in a $100 \times 100$ grid. As for the potential LM locations, $5 \times |V|$ locations on the 100$ \times $100 grid were chosen at random. In total, we had $180$ instances on which all the experiments were performed. The sensing range for the bearing sensors, $\rho_s$, was fixed at $35$ units.

\subsection{Branch-and-cut algorithm performance}

The branch-and-cut algorithm presented in this article was implemented in the Julia programming language using the callback functionality of CPLEX. The internal cut-generation routines of CPLEX were switched off and CPLEX was used only to manage the enumeration tree in the branch-and-cut algorithm. All computation times are reported in seconds. The performance of the algorithm was tested on randomly generated test instances. The branch-and-cut algorithm, to compute optimal solution to the problem, is very effective in computing optimal solutions for all the instances that were generated, within than 30 targets. The computation time for most of the instances was less than two seconds and at most five seconds; hence, for all practical routing purposes the algorithm can be deployed online for $ |V|\leqslant 30$.

\begin{table}[]
    \centering
    \begin{tabular}{c|c|c|c}
         \toprule 
         \# vehicles & \# $|V|$ & \# constraints & \# landmarks  \\
         \midrule 
         2 & 20 & 19.8 & 9.4 \\
         2 & 25 & 24.6 & 9 \\
         2 & 30 & 138.2 & 8.8 \\
         3 & 20 & 25.2 & 9.2 \\
         3 & 25 & 56.7 & 8.6 \\
         3 & 30 & 96.8 & 8.9 \\
         4 & 20 & 17.6 & 9 \\
         4 & 25 & 52.4 & 8.5 \\
         4 & 30 & 82 & 9.1 \\
         \bottomrule
    \end{tabular}
    \caption{Average number of sub-tour and path elimination constraints in Eq. \eqref{e:2}--\eqref{e:4} and the average number of landmarks placed in the optimal solution to the MVPLC instances.}
    \label{tab:1}
\end{table}

The table \ref{tab:1} shows the average number of ``user defined cuts'' -- sub-tour elimination constraints in Eq. \eqref{e:2} and path elimination constraints in Eqs. \eqref{e:3} and \eqref{e:4} added in S\textsc{tep} 4 of the branch-and-cut algorithm detailed in Sec. \ref{sec:bandc}, and the average number of LMs required in the optimal solution; this indicates the effectiveness of the branch-and-cut approach. Several conclusions can be drawn from the results in table \ref{tab:1}.

\begin{enumerate}
\item When fixing the number of UVs and increasing the number of targets, the number of cuts, especially sub-tour elimination constraints in \eqref{e:2}, will increase rapidly, which is consistent with typical multiple traveling salesmen problems.
\item Almost regardless of the number of UVs and targets, the average number of LMs remains fairly steady, indicating that it might be only a function of the area of the grid, and the sensing range, which were both fixed in our setup.
\item It is worth recalling that in the mathematical formulation solved by our algorithm, the cost for placing a landmark was set to be unity, which is way smaller than the cost for traveling through almost any edge. Even in this setting, the optimal solution to every randomly generated test instance only included placing a few landmarks. Thus, it is almost safe to conclude that the MVPLC we are studying is primarily an optimal routing problem, and its solution favors a minimal set of locations for placing the necessary landmarks. This need not necessarily be the case when the landmark placement costs are comparable to the travel costs. 
\end{enumerate}

\subsection{Simulation results}

For the simulation \emph{i.e.}, online estimation algorithm using the results of the branch-and-cut algorithm, we consider one instance with 15 waypoints (targets) and two vehicles. The routes for each vehicle and the locations where landmarks are to be placed are obtained using the branch-and-cut algorithm. For the simulation, we assume that the vehicles have a very high turn rate such that on reaching a way point, they can immediately point towards the next way point. Furthermore, $\rho_s$ was set to $35$ units. The estimated states for the vehicles were used in the way point controller instead of the true states to show that the vehicles can indeed travel in a GPS-restricted environment provided that the condition for path to at least two LMs is always maintained. The simulations were ran for $3000$ iterations. Without loss of generality, for the purpose of simulation, the unit for distance and time were chosen as meters (m) and seconds (s), respectively.

\begin{figure}[h!]
\includegraphics[width=80mm,scale=0.4]{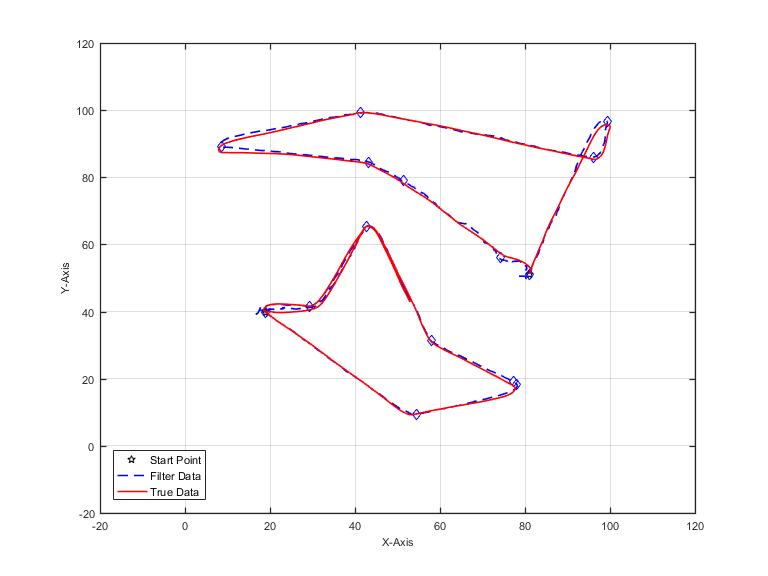}
\includegraphics[width=80mm,scale=0.4]{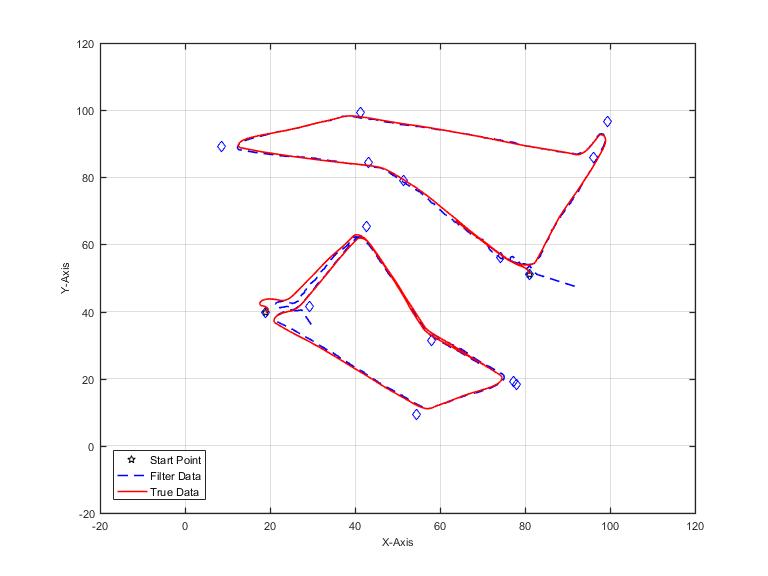}
\caption{Plot showing trajectories of the two vehicles with sensing range as 35 units and for $|V|$ = 15. The first plot shows the trajectories for a controller gain of $2.0$ and the second plot shows the trajectories for a controller gain of $0.4$}
\label{fig:simulation_1}
\end{figure}

The plots of the true and estimated trajectories for the instances with two different values for controller gains are shown in the Fig. \ref{fig:simulation_1}.

\begin{figure}[h!]
\includegraphics[width=80mm,scale=0.4]{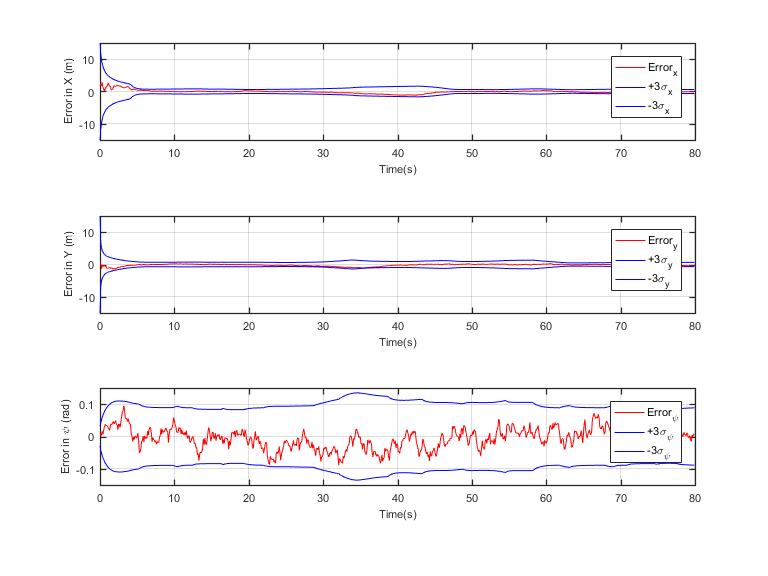}
\includegraphics[width=80mm,scale=0.4]{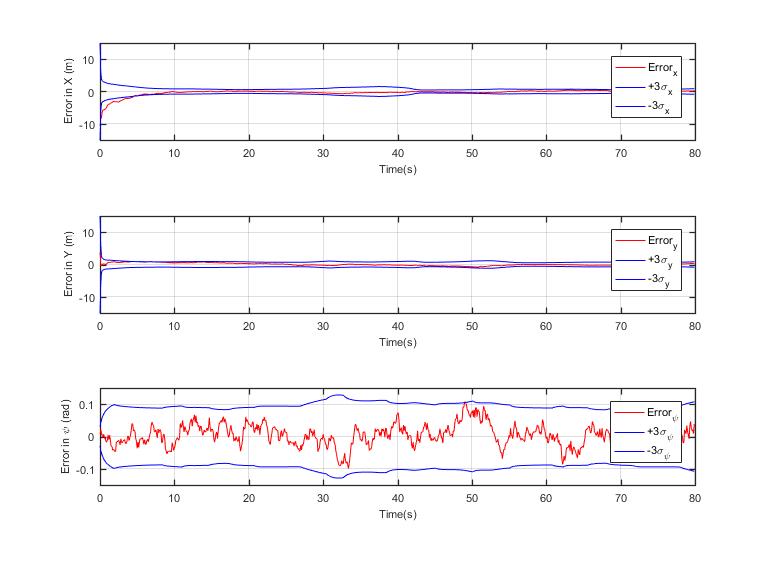}
\caption{Plot showing the error in X direction, Y direction and heading ($\psi$) along with their respective $3\sigma$ bounds for the two chosen values of controller gains with $|V|$ = 15.}
\label{fig:simulation_2}
\end{figure}

The next set of plots in Fig. \ref{fig:simulation_2} shows the actual errors (and the corresponding $3\sigma$ bounds) in position and heading for different values of controller gains.

The above simulation results illustrate that the mathematical formulation presented in this article is able to capture the localization constraints in an effective manner.

\subsection{Experimental Results}
To further validate our methodology and solution, we performed field experiments for one test instance, routing two Turtlebots. The update step for the EKF requires the exteroceptive sensor mounted on a vehicle to detect/communicate with external environment agents frequently. In this case, the external environment agents are simulated Landmarks (LMs) and the exteroceptive sensor is an Orbbec Astra camera mounted on the Turtlebot 2.0 platform. The Orbbec Astra camera is a stereo camera, which can measure depth alongside providing an RGB stream. However, a bearing-only sensor is used for the update step, thus, obviating the use of stereo capability of the camera. The camera has a $60 ^{\circ} $ horizontal field of view with an image sensing range of $0.6$m to $8$m. It weighs $0.3$ kg with a dimension of $165 \times 30 \times 40$ mm (with Turtlebot 2.0 mounting capability). It also provides an RGB stream of size $640 \times 480$ pixels at $30$ frames per second, making it ideal for our use.

A testing field of size $4 \times 8$ m is utilized. For proper detection, the sensing range regarding the simulated landmarks is no greater than $2$m. When the sensing range is between $1.2$m and $2$m, it is observed that the probability of detection gets around $0.6$. However, due to the nature of real-world vehicle movement, if the detection of landmarks fails at some point, the vehicle might deviate quite much from the desirable trajectory, such that it might miss the next waypoint; or, in extreme cases, it might run into some area that no landmark could be detected at all, and fail the entire mission. To avoid taking risk, while still utilizing the optimization model, a considerable safe margin is introduced, that the sensing range is set to be $1.2$m when solving the model; in this way, any landmark placed will be detectable.

\begin{figure}
\includegraphics[width=80mm,scale=0.4]{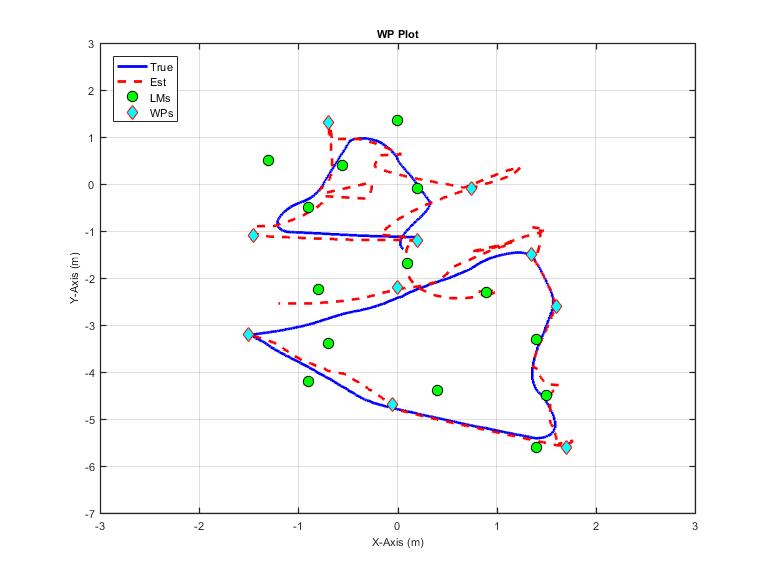}
\includegraphics[width=80mm,scale=0.4]{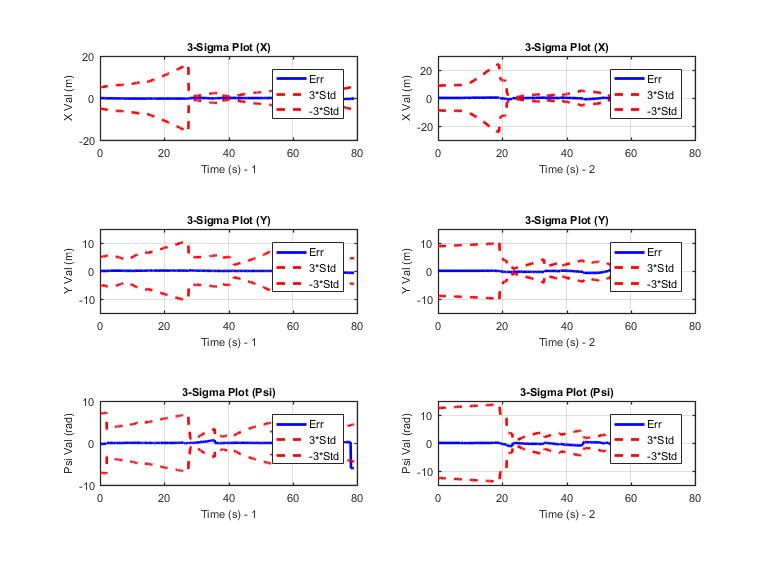}
\caption{Plot showing actual trajectories, and the error in X direction, Y direction and heading ($\psi$) along with their respective $3\sigma$ bounds of the two Turtlebots, routing using estimated states.}
\label{fig:exp}
\end{figure}

For the hardware experiments, a high controller gain was chosen to ensure a high turn rate for the vehicles. The vehicles' velocity was kept constant at $0.2$ m/s. A $p$ value of $0.05$ m (collision avoidance distance) was used for the separation distance between the vehicles and any location where the landmarks were to be placed. Although this distance can be increased, a small value was chosen to make the instances feasible. Using the landmarks at the locations obtained from the algorithm, an EIF was used to estimate the position and heading of the vehicles using the bearing measurements. The results are shown in Fig. \ref{fig:exp}. It can be seen that the error lies within the $3\sigma$ bound most of the time, hence ensuring that the placement of landmarks obtained from the optimization model is sufficient for localization of the vehicles.

The substantial difference between the estimated and the true trajectories can be attributed to the fact that the formulation presented in this article considers a $360^{\circ}$ field of view of the camera while in reality, the field of view is around $60^{\circ}$. Nevertheless, the errors are within the $3\sigma$ bounds even for such an approximation.

\section{Conclusion}
In this article, a systematic method to address the problem of joint routing and localization for multiple UVs in a GPS-denied or GPS-restricted environment was presented. The optimization problem computes routes and determines the minimal set of locations where landmarks need to be placed to enable vehicle localization. This solution is embedded with estimation algorithms to estimate the states of the vehicles. Fast separation and branch-and-cut algorithms to compute an optimal solution have also been developed. The proposed system architecture have been tested extensively via simulation and field experiments, and the feasibility has clearly been validated. Future work can be focussed on handling field of view constraints for the camera in the framework presented in this article and taking into account more realistic constraints for the fleet of the vehicles like fuel constraints or heterogeneous constraints applied on specific vehicles.

\bibliography{bibtex_database2}
\bibliographystyle{plain}

\end{document}